\documentclass[reqno]{amsart}

\usepackage{lmodern}
\usepackage[T1]{fontenc}
\usepackage[utf8]{inputenc}
\usepackage[english]{babel}
\usepackage{microtype} 

\usepackage{mathdots,latexsym,mathtools,amssymb,booktabs,tikz}

\usepackage[capitalize,noabbrev]{cleveref}
\usepackage[centertableaux]{ytableau}

\newtheorem{theorem}{Theorem}
\newtheorem{proposition}[theorem]{Proposition}

\newtheorem{conjecture}[theorem]{Conjecture}

\theoremstyle{definition}

\newtheorem{example}[theorem]{Example}

\newcommand{\symS}{\mathfrak{S}}

\newcommand{\setZ}{\mathbb{Z}}
\newcommand{\setR}{\mathbb{R}}
\newcommand{\setQ}{\mathbb{Q}}
\newcommand{\setC}{\mathbb{C}}

\newcommand{\GT}{\mathcal{GT}}
\newcommand{\SSYT}{\mathrm{SSYT}}
\newcommand{\schurS}{\mathrm{s}}
\newcommand{\key}{\kappa}
\newcommand{\ehr}{i}
\newcommand{\type}{\mathrm{type}}

\newcommand{\defin}[1]{\emph{#1}}

\title{Ehrhart positivity and Demazure characters}

\author{Per Alexandersson}
\address{Department of Mathematics, KTH,\\
SE-100 44 Stockholm, Sweden}
\email{per.w.alexandersson@gmail.com}

\author{Elie Alhajjar}
\address{Department of Mathematics, USMA,\\
West Point, New York, USA}
\email{eliealhajjar@gmail.com}

\begin{document}

\maketitle

\begin{abstract}

Demazure characters, also known as key polynomials, generalize the classical Schur polynomials. In particular, when all variables are set equal to $1$, these polynomials count the number of integer points in a certain class of Gelfand--Tsetlin polytopes. This property highlights the interaction between the corresponding polyhedral and combinatorial structures via Ehrhart theory. In this paper, we give an overview of results concerning the interplay between the geometry of Gelfand--Tsetlin polytopes and their Ehrhart polynomials. Motivated by strong computer evidence, we propose several conjectures about the non-negativity of the coefficients of such polynomials.

\end{abstract}



\section{Introduction}\label{sec:intro}

The theory of Schur polynomials can be seen from two different sets of lenses. 
On one hand, the traditional approach begins with a definition involving quotients of matrix determinants. 
This method is mainly useful in representation theory, since it is derived as a special case of the 
Weyl character formula. On the other hand, the combinatorial approach uses the sum expansion over semi-standard 
Young tableaux of fixed shape. Note that it is not hard to show explicitly the equivalence of these two approaches.

Demazure characters \cite{Demazure1974nouvelle,Demazure1974}, also known as key polynomials, generalize the classical Schur polynomials.
Key polynomials can be computed recursively via divided difference operators,
and are closely related to Schubert polynomials. In particular, each Schubert polynomial can be expressed 
as a non-negative integer combination of key polynomials.

A combinatorial formula using semi-standard Young tableaux was discovered in \cite{Lascoux1990Keys},
and this is where the notion of key polynomials come from.
Key polynomials are specializations of non-symmetric Macdonald polynomials, \cite{HaglundNonSymmetricMacdonald2008}
so the combinatorial formula of J.~Haglund gives an alternative formulation, using \emph{skyline fillings}.
S.~Mason\cite{Mason2009} explores two variations of skyline fillings, both of which give key polynomials. 
It is possible to interpolate between the two skyline models, see \cite{Kurland2016,AlexanderssonSawhney18}.

For each a partition $\lambda$, one can construct a Gelfand--Tsetlin polytope.
These polytopes play a crucial role in representation theory, algebraic geometry and combinatorics. 
Their importance stems from the fact that their integer points are in bijection with semi-standard Young tableaux.

R.~Stanley and A.~Postnikov\cite{PostnikovStanley2008} study a certain subfamily of key polynomials,
and prove that these are flagged Schur polynomials. This implies that key polynomials in this family 
can be computed as a certain sum over lattice points in a single face of a Gelfand--Tsetlin polytope.

The result by R.~Stanley and A.~Postnikov can be extended to the full family of key polynomials;
V.~Kirichenko, E.~Smirnov and V.~Timorin \cite{Kiritchenko2010} show that
key polynomials can be expressed as a sum over lattice pooints in a certain union of faces in a Gelfand--Tsetlin polytope.
This way of thinking about key polynomials is not as well-known as the other interpretations,
and the purpose of this survey is to emphasize the polyhedral aspect of key polynomials.
A related result appears in \cite{FeiginMakhlin2016}, where Hall--Littlewood polynomials 
are expressed as a weighted sum over lattice points in Gelfand--Tsetlin polytopes.

Recent research has been focused on products involving key polynomials, 
see \cite{Pun2016Thesis}. The main motivation for studying key polynomials is to gain insight about the expression of products of Schubert polynomials in terms of Schubert polynomials, a problem of main importance in representation theory. The close relationship between key and Schubert polynomials is emphasized in \cite{ReinerShimozono1995}.

The purpose of the current paper is two-fold$\colon$on one hand, we aim to collect some of the main results related to the study of key polynomials. On the other hand, we propose several conjectures concerning the non-negativity of the coefficients of the `stretched' version of such polynomials. 
In \cref{sec:prelim} below, we provide the basic material and fix the terminology for the remainder of the paper. 
\cref{sec:reducedKogan} deals with the facial description of GT-polytopes and the formal definition of key polynomials, where we give several examples to illustrate the main ideas. 
In \cref{sec:ehrhart} we introduce the connection to Ehrhart theory through Kostka coefficients and in 
\cref{sec:conjectures} we mention a sample of the computations that lead eventually to the main conjecture. 
Finally, we reconstruct a counterexample in \cref{sec:faces} that shows the failure of the non-negativity argument in the case of arbitrary faces of GT-polytopes.

\section{Preliminaries}\label{sec:prelim}

Given an integer partition $\lambda_1 \geq \lambda_2 \geq \lambda_3\geq \dotsb$,
we can associate a \emph{Young diagram} of shape $\lambda$ as a diagram in the plane with $\lambda_i$ left-justified boxes in row $i$.
For example, $\lambda=(5,3,2,2)$ gives rise to the following Young diagram:
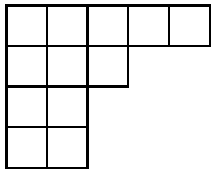
\[
 \begin{ytableau}
  \; & \; & \; & \; & \; \\
  \; & \; & \; \\
  \; & \; \\
  \; & \; \\
 \end{ytableau}
\]
A \emph{semi-standard Young tableau} of shape $\lambda$ is an assignment of natural numbers to the boxes of the Young diagram,
such that rows are weakly increasing from left to right, and columns are strictly increasing from top to bottom.
Only the first of the following three assignments is a semi-standard Young tableau
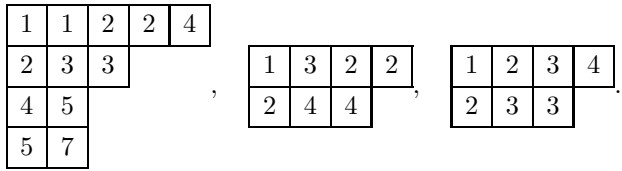
\[
 \begin{ytableau}
  1 & 1 & 2 & 2 & 4 \\
  2 & 3 & 3 \\
  4 & 5 \\
  5 & 7 \\
 \end{ytableau}
 ,
 \quad
  \begin{ytableau}
  1 & 3 & 2 & 2 \\
  2 & 4 & 4 \\
 \end{ytableau}
 ,
 \quad
  \begin{ytableau}
  1 & 2 & 3 & 4 \\
  2 & 3 & 3 \\
 \end{ytableau}.
\]

\subsection{Gelfand--Tsetlin polytopes}\label{ssec:gtPolytope}

There are several families of polytopes which are referred to as Gelfand--Tsetlin polytopes, see for example
\cite{GusevKiritchenkoTimorin2013} and \cite{DeLoeraMcAllister2004,King04stretched}.
A \emph{Gelfand--Tsetlin pattern} or GT-pattern for short is a triangular array $(x_{ij})$
visualized as
\begin{equation}\label{eq:gtTriangular}
\scriptstyle \begin{array}{cccccccccccccc}
x_{n1} & & x_{n2} & & \cdots  & & x_{nn} \\
& \ddots & & \ddots &  & \iddots &   \\
 &  &   x_{21} &  & x_{22}  \\
  &  &    & x_{11}
\end{array}
\end{equation}
satisfying the inequalities
\begin{equation} \label{eq:gtinequalities}
x_{i+1,j} \geq x_{ij} \text{ and } x_{ij} \geq x_{i+1,j+1}
\end{equation}
for all values of $i$, $j$ where the indexing is defined.
The inequalities simply state that down-right diagonals 
are weakly decreasing and down-left diagonals are weakly increasing.

Given an integer partition $\lambda$, the \emph{Gelfand--Tsetlin polytope} $\GT(\lambda) \subset \setR^{\tfrac{n(n+1)}{2}}$ is 
the convex polytope of Gelfand--Tsetlin patterns defined by the inequalities in \cref{eq:gtinequalities}
together with the equalities $x_{ni} = \lambda_i$ for $i=1,2,\dots,n$.


The polytope $\GT(\lambda)$ has integer vertices. 
In fact, it has a unimodular triangulation, see \cite{Alexandersson16Counterexamples}.
Also, note that $k\cdot \GT(\lambda) = \GT(k\lambda)$ for all $k\geq 0$.

\subsection{Bijection with semi-standard Young tableaux}\label{subsec:bijection}

Note that \eqref{eq:gtinequalities} implies that any two adjacent rows in an integral GT-pattern form a skew Young diagram,
see the standard textbook by R.~Stanley\cite{StanleyEC2} for terminology.

This property enables us to define a bijection with Young tableaux ---
the skew shape defined by row $j$ and $j+1$ in an integral GT-pattern $G$
describes which boxes in a Young tableau $T$ have content $j$.
In particular, tableaux of shape $\lambda$ are in bijection
with integral GT-patterns with topmost row equal to $\lambda$.
See \cref{fig:gtexample} for an example of this correspondence.
\begin{figure}[!ht]
\centering
\setcounter{MaxMatrixCols}{20}
\begin{equation}\label{eq:triangulargt}
\begin{matrix}
5 &   & 4 &   & 2 &   & 1 &   & 1 &  & 0\\
 & 5 &   & 3 &   & 2 &   & 1 &   & 0\\
 &  & 3 &   & 3 &   & 2 &   & 1  \\
 &  &  & 3 &   & 3 &   & 1 &  \\
 &  &  &  & 3 &   & 2 &  \\
 &  &  &  &  & 3 & \\
\end{matrix}
\quad
\longleftrightarrow
\quad
\begin{ytableau}
1&1&1&5&5 \\
2&2&3&6\\
3&4\\
4\\
6
\end{ytableau}
\end{equation}
\caption{The GT-pattern corresponding to a Young tableau.
For example, the third row tells us that the shape of the entries $\leq 3$ in the tableau is $(3,3,1)$.
}\label{fig:gtexample}
\end{figure}

Note that in any integral GT-pattern, $x_{i+1,j} - x_{ij}$ counts the number of
boxes with content $i$ in row $j$ in the corresponding tableau.
Given a GT-pattern $G$, we define the \defin{weight} $w(G)$ as the vector
\begin{equation}\label{eq:weight}
 w_i(G) \coloneqq \sum_{j=1}^{i+1} x_{i+1,j} - \sum_{j=1}^{i} x_{ij},
\end{equation}
where $x_{0j}\coloneqq 0$.
Thus, an integral GT-pattern with weight $w$ is in bijection with a semi-standard Young tableau 
with $w_i$ entries equal to $i$. 
Hence, integer points in $\GT(\lambda)$ are in bijection with $\SSYT(\lambda,n)$ --- the set of 
semi-standard Young tableaux of shape $\lambda$ with maximal entry $n$.
Given $\lambda$ and $w$, let $\GT(\lambda,w) \subseteq \GT(\lambda)$ be the intersection 
of $\GT(\lambda)$ with the hyperplanes defined by \eqref{eq:weight}.
The lattice points in $\GT(\lambda,w)$ are enumerated by the Kostka coefficients, see \cref{ssec:kostka}
below.

One can then define the \defin{Schur polynomials} $\schurS_\lambda(z_1,\dotsc,z_n)$ as 
\begin{align}\label{eq:schurGT}
 \schurS_\lambda(z_1,\dotsc,z_n) \coloneqq 
 \sum_{G \in \GT(\lambda) \cap \setZ^{\tfrac{n(n+1)}{2}} } z_1^{w_1(G)} \dotsm z_n^{w_n(G)}.
\end{align}

In particular, by using the Weyl dimension formula \cite[Eq. 7.105]{StanleyEC2} we have that
\begin{equation}\label{eq:weylDimFormula}
 \schurS_\lambda(\underbrace{1,1,\dotsc,1}_n) = |\GT(\lambda) \cap \setZ^{\tfrac{n(n+1)}{2}}| = 
 |\SSYT(\lambda,n)|
 = \prod_{1 \leq i < j \leq n} \frac{\lambda_i - \lambda_j + j-i}{j-i}.
\end{equation}

The fact that the right hand side is a polynomial in the $\lambda_i$ is not 
a complete surprise. Gelfand--Tsetlin polytopes belong to a larger family of so called \emph{marked order polytopes},
generalizing the notion of order polytopes introduced by R.~Stanley in \cite{Stanley86TwoPosetPolytopes}.
K.~Jochemko and R.~Sanyal proved in \cite{JochemkoSanyal2014} that marked order polytopes give rise to Ehrhart 
functions which are piecewise polynomial in the markings, which in our case are the parts of $\lambda$.

\subsection{Skew GT-polytopes}

In analogy with the triangular case, one can define 
so called skew Gelfand--Tsetlin polytopes using parallelogram arrangements of non-negative numbers,
\begin{equation}\label{eq:gtPaarallelogram}
\scriptstyle \begin{array}{cccccccccccccc}
x_{n1} & & x_{n2} & & \cdots  & & \cdots  & & x_{nm} \\
& \ddots & & \ddots &  & \ddots && \ddots && \ddots &   \\
 &  &   x_{11} &  & x_{12}& & \cdots  & & \cdots  & & x_{1m}  \\
  &  &    & x_{01} &  & x_{02}& & \cdots  & & \cdots  & & x_{0m} 
\end{array}
\end{equation}
satisfying the inequalities
\begin{equation} \label{eq:gtinequalitiesP}
x_{i+1,j} \geq x_{ij} \text{ and } x_{ij} \geq x_{i+1,j+1}
\end{equation}
for all $i$, $j$ where the indexing is defined.

Consider an $(n+1)\times m$ GT-pattern with top row $\lambda$ and bottom row $\mu$ --- that is, $x_{ni}=\lambda_i$ and $x_{0i}=\mu_i$ for $i=1,\dotsc,m$. 
These equalities together with the above inequalities define a convex polytope, 
the \emph{skew Gelfand--Tsetlin polytope}, $\GT(\lambda/ \mu) \subset \mathbb{R}^{(n+1)m}$. 
Note that the vertices of such polytopes have integer coordinates. In fact, the skew Gelfand--Tsetlin polytopes admit 
a unimodular triangulation.

The weight of a GT-pattern in the parallelogram form is defined in the same way as in the triangular form \eqref{eq:weight}. 
Similar to the bijection described in \cref{subsec:bijection}, the integer points in $\GT(\lambda/ \mu)$ correspond to 
skew Young tableaux with shape $\lambda/\mu$, where the entries belong to the set $\{1, 2,\dots, n\}$.
This allow us to define the \emph{skew Schur polynomials} as 
\begin{align}\label{eq:skewschurGT}
 \schurS_{\lambda/\mu}(z_1,\dotsc,z_n) \coloneqq 
 \sum_{G \in \GT(\lambda/\mu) \cap \setZ^{(n+1)m} } z_1^{w_1(G)} \dotsm z_n^{w_n(G)}.
\end{align}
We note that there is no simple formula for computing the specialization $\schurS_{\lambda/\mu}(1,1,\dotsc,1)$ in general.

\medskip 

By specifying a weight vector $w$, we intersect $\GT(\lambda/ \mu)$ with a set of hyperplanes,
and denote the resulting polytope by $\GT(\lambda/ \mu, w)$. 
As with $\GT(\lambda,w)$, the polytopes $\GT(\lambda/ \mu, w)$ are also not integral in general.

\section{Reduced Kogan faces and key polynomials}\label{sec:reducedKogan}

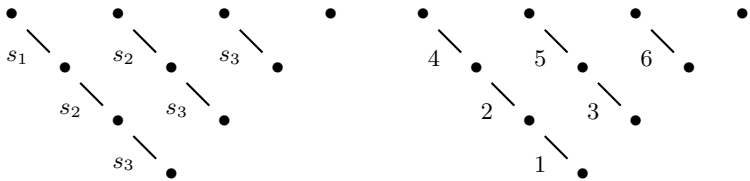
\begin{figure}[!ht]
\centering
\begin{tikzpicture}
\node (N11) at (0,0) {$\bullet$};
\node [below right of=N11] (N21) {$\bullet$};
\node [above right of=N21] (N12) {$\bullet$};
\node [below right of=N12] (N22) {$\bullet$};
\node [above right of=N22] (N13) {$\bullet$};
\node [below right of=N13] (N23) {$\bullet$};
\node [above right of=N23] (N14) {$\bullet$};
\node [below right of=N21] (N31) {$\bullet$};
\node [below right of=N22] (N32) {$\bullet$};
\node [below right of=N31] (N41) {$\bullet$};

\path[-,thick,every node/.style={font=\sffamily\small}]
(N11) edge node [below left] {$s_1$} (N21)
(N21) edge node [below left] {$s_2$} (N31)
(N31) edge node [below left] {$s_3$} (N41)
(N12) edge node [below left] {$s_2$} (N22)
(N22) edge node [below left] {$s_3$} (N32)
(N13) edge node [below left] {$s_3$} (N23);
\end{tikzpicture}
\qquad
\begin{tikzpicture}
\node (N11) at (0,0) {$\bullet$};
\node [below right of=N11] (N21) {$\bullet$};
\node [above right of=N21] (N12) {$\bullet$};
\node [below right of=N12] (N22) {$\bullet$};
\node [above right of=N22] (N13) {$\bullet$};
\node [below right of=N13] (N23) {$\bullet$};
\node [above right of=N23] (N14) {$\bullet$};
\node [below right of=N21] (N31) {$\bullet$};
\node [below right of=N22] (N32) {$\bullet$};
\node [below right of=N31] (N41) {$\bullet$};

\path[-,thick,every node/.style={font=\sffamily\small}]
(N11) edge node [below left] {$4$} (N21)
(N21) edge node [below left] {$2$} (N31)
(N31) edge node [below left] {$1$} (N41)
(N12) edge node [below left] {$5$} (N22)
(N22) edge node [below left] {$3$} (N32)
(N13) edge node [below left] {$6$} (N23);
\end{tikzpicture}
\caption{
The transposition corresponding to different equalities between
entries in a Kogan face and the reading order of these.
}\label{fig:reducedKoganFace}
\end{figure}

By imposing some additional equalities on the coordinates of a GT-polytope, 
one can obtain faces of the polytope.
There is a particular interest with equalities of the form $x_{ij}=x_{i,j+1}$.
By imposing a set of such equalities, we obtain a \defin{Kogan face} of the GT-polytope.
To each equality of the form $x_{ij}=x_{i,j+1}$, we associate the transposition $s_{n-i+j-1}$,
as shown in \cref{fig:reducedKoganFace}.
We then construct a word from these transpositions by reading the entries 
from bottom to top, left to right.
If this word is \emph{reduced}, we say that the corresponding Kogan face is reduced.
Note that the same word might be constructed from equalities in several different ways.
The \defin{type of a Kogan face} is the permutation obtained from the word.
We should really view the equalities that define a Kogan face as some special set of equalities ---
a point in this face might satisfy some additional equalities present, if it is also a member of some sub-face.
This implies that a point in the GT-polytope can be a member of \emph{several} (reduced) Kogan faces.
For example, the point where all equalities are present is the unique face of type $w_0$, the longest permutation in $\symS_n$.
This point is a sub-face of all other Kogan faces.
Furthermore, the full GT-polytope has the identity permutation as type.

\begin{example}
Consider the face in \cref{fig:reducedKoganFaceExample}.
All marked equalities are in the south-east direction, so we obtain the word $s_3s_1s_2s_3$.
It is straightforward to verify that this word is reduced, so this face is a reduced Kogan face.
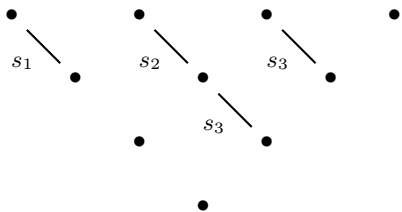
\begin{figure}[!ht]
\centering
\begin{tikzpicture}[node distance =1.2cm]
\node (N11) at (0,0) {$\bullet$};
\node [below right of=N11] (N21) {$\bullet$};
\node [above right of=N21] (N12) {$\bullet$};
\node [below right of=N12] (N22) {$\bullet$};
\node [above right of=N22] (N13) {$\bullet$};
\node [below right of=N13] (N23) {$\bullet$};
\node [above right of=N23] (N14) {$\bullet$};
\node [below right of=N21] (N31) {$\bullet$};
\node [below right of=N22] (N32) {$\bullet$};
\node [below right of=N31] (N41) {$\bullet$};

\path[-,thick,every node/.style={font=\sffamily\small}]
(N11) edge node [below left] {$s_1$} (N21)
(N12) edge node [below left] {$s_2$} (N22)
(N22) edge node [below left] {$s_3$} (N32)
(N13) edge node [below left] {$s_3$} (N23);
\end{tikzpicture}
\caption{
A reduced Kogan face with the word $\omega = s_3s_1s_2s_3$.
}\label{fig:reducedKoganFaceExample}
\end{figure}

Suppose we wish to examine the GT-polytope $\GT(\lambda)$ with $\lambda=(4,3,3,2)$.
The lattice points in the reduced Kogan face with the word $s_3s_1s_2s_3$ are the following GT-patterns:
\begin{equation*}
\begin{tikzpicture}[scale=0.4]
\node at (3,0) {$4$};
\node at (5,0) {$3$};
\node at (7,0) {$3$};
\node at (9,0) {$2$};
\node at (4,-1) {$4$};
\node at (6,-1) {$3$};
\node at (8,-1) {$3$};
\node at (5,-2) {$3$};
\node at (7,-2) {$3$};
\node at (6,-3) {$3$};
\end{tikzpicture}
\qquad
\begin{tikzpicture}[scale=0.4]
\node at (3,0) {$4$};
\node at (5,0) {$3$};
\node at (7,0) {$3$};
\node at (9,0) {$2$};
\node at (4,-1) {$4$};
\node at (6,-1) {$3$};
\node at (8,-1) {$3$};
\node at (5,-2) {$4$};
\node at (7,-2) {$3$};
\node at (6,-3) {$3$};
\end{tikzpicture}
\qquad
\begin{tikzpicture}[scale=0.4]
\node at (3,0) {$4$};
\node at (5,0) {$3$};
\node at (7,0) {$3$};
\node at (9,0) {$2$};
\node at (4,-1) {$4$};
\node at (6,-1) {$3$};
\node at (8,-1) {$3$};
\node at (5,-2) {$4$};
\node at (7,-2) {$3$};
\node at (6,-3) {$4$};
\end{tikzpicture}
\end{equation*}
\end{example}

It is evident from the bubble-sort algorithm \cite{Knuth1998ArtOfProgramming} that 
there is at least one reduced Kogan face of every type $\sigma \in \symS_n$.

\begin{proposition}[See \cite{Kogan2000schubertgeometry}]\label{prop:kempfCase}
If $\sigma \in \symS_n$ avoids the permutation pattern $132$ (such permutations are known as Kempf permutations),
then there is a unique reduced Kogan face of type $\sigma$.
\end{proposition}
See \cref{ssec:determinant} for more background on this special case.

\begin{example}
There are $11$ reduced Kogan faces for $n=4$ that are covered by \cref{prop:kempfCase}.
These faces are illustrated in \cref{fig:kempf}.
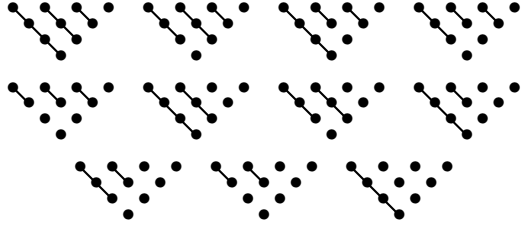
\begin{figure}[!ht]
\centering

\begin{tikzpicture}[node distance=0.3cm]
\node (N11) at (0,0) {$\bullet$};
\node [below right of=N11] (N21) {$\bullet$};
\node [above right of=N21] (N12) {$\bullet$};
\node [below right of=N12] (N22) {$\bullet$};
\node [above right of=N22] (N13) {$\bullet$};
\node [below right of=N13] (N23) {$\bullet$};
\node [above right of=N23] (N14) {$\bullet$};
\node [below right of=N21] (N31) {$\bullet$};
\node [below right of=N22] (N32) {$\bullet$};
\node [below right of=N31] (N41) {$\bullet$};
\path[-,thick,every node/.style={font=\sffamily\small}]
(N11) edge node [below left] {} (N21)
(N12) edge node [below left] {} (N22)
(N13) edge node [below left] {} (N23)
(N21) edge node [below left] {} (N31)
(N22) edge node [below left] {} (N32)
(N31) edge node [below left] {} (N41)
;
\end{tikzpicture}
%
\begin{tikzpicture}[node distance=0.3cm]
\node (N11) at (0,0) {$\bullet$};
\node [below right of=N11] (N21) {$\bullet$};
\node [above right of=N21] (N12) {$\bullet$};
\node [below right of=N12] (N22) {$\bullet$};
\node [above right of=N22] (N13) {$\bullet$};
\node [below right of=N13] (N23) {$\bullet$};
\node [above right of=N23] (N14) {$\bullet$};
\node [below right of=N21] (N31) {$\bullet$};
\node [below right of=N22] (N32) {$\bullet$};
\node [below right of=N31] (N41) {$\bullet$};
\path[-,thick,every node/.style={font=\sffamily\small}]
(N11) edge node [below left] {} (N21)
(N12) edge node [below left] {} (N22)
(N13) edge node [below left] {} (N23)
(N21) edge node [below left] {} (N31)
(N22) edge node [below left] {} (N32)
%
;
\end{tikzpicture}
%
\begin{tikzpicture}[node distance=0.3cm]
\node (N11) at (0,0) {$\bullet$};
\node [below right of=N11] (N21) {$\bullet$};
\node [above right of=N21] (N12) {$\bullet$};
\node [below right of=N12] (N22) {$\bullet$};
\node [above right of=N22] (N13) {$\bullet$};
\node [below right of=N13] (N23) {$\bullet$};
\node [above right of=N23] (N14) {$\bullet$};
\node [below right of=N21] (N31) {$\bullet$};
\node [below right of=N22] (N32) {$\bullet$};
\node [below right of=N31] (N41) {$\bullet$};
\path[-,thick,every node/.style={font=\sffamily\small}]
(N11) edge node [below left] {} (N21)
(N12) edge node [below left] {} (N22)
(N13) edge node [below left] {} (N23)
(N21) edge node [below left] {} (N31)
%
(N31) edge node [below left] {} (N41)
;
\end{tikzpicture}
%
\begin{tikzpicture}[node distance=0.3cm]
\node (N11) at (0,0) {$\bullet$};
\node [below right of=N11] (N21) {$\bullet$};
\node [above right of=N21] (N12) {$\bullet$};
\node [below right of=N12] (N22) {$\bullet$};
\node [above right of=N22] (N13) {$\bullet$};
\node [below right of=N13] (N23) {$\bullet$};
\node [above right of=N23] (N14) {$\bullet$};
\node [below right of=N21] (N31) {$\bullet$};
\node [below right of=N22] (N32) {$\bullet$};
\node [below right of=N31] (N41) {$\bullet$};
\path[-,thick,every node/.style={font=\sffamily\small}]
(N11) edge node [below left] {} (N21)
(N12) edge node [below left] {} (N22)
(N13) edge node [below left] {} (N23)
(N21) edge node [below left] {} (N31)
%
;
\end{tikzpicture}

%
\begin{tikzpicture}[node distance=0.3cm]
\node (N11) at (0,0) {$\bullet$};
\node [below right of=N11] (N21) {$\bullet$};
\node [above right of=N21] (N12) {$\bullet$};
\node [below right of=N12] (N22) {$\bullet$};
\node [above right of=N22] (N13) {$\bullet$};
\node [below right of=N13] (N23) {$\bullet$};
\node [above right of=N23] (N14) {$\bullet$};
\node [below right of=N21] (N31) {$\bullet$};
\node [below right of=N22] (N32) {$\bullet$};
\node [below right of=N31] (N41) {$\bullet$};
\path[-,thick,every node/.style={font=\sffamily\small}]
(N11) edge node [below left] {} (N21)
(N12) edge node [below left] {} (N22)
(N13) edge node [below left] {} (N23)
%
%
;
\end{tikzpicture} 
%
\begin{tikzpicture}[node distance=0.3cm]
\node (N11) at (0,0) {$\bullet$};
\node [below right of=N11] (N21) {$\bullet$};
\node [above right of=N21] (N12) {$\bullet$};
\node [below right of=N12] (N22) {$\bullet$};
\node [above right of=N22] (N13) {$\bullet$};
\node [below right of=N13] (N23) {$\bullet$};
\node [above right of=N23] (N14) {$\bullet$};
\node [below right of=N21] (N31) {$\bullet$};
\node [below right of=N22] (N32) {$\bullet$};
\node [below right of=N31] (N41) {$\bullet$};
\path[-,thick,every node/.style={font=\sffamily\small}]
(N11) edge node [below left] {} (N21)
(N12) edge node [below left] {} (N22)
%
(N21) edge node [below left] {} (N31)
(N22) edge node [below left] {} (N32)
(N31) edge node [below left] {} (N41)
;
\end{tikzpicture} 
%
\begin{tikzpicture}[node distance=0.3cm]
\node (N11) at (0,0) {$\bullet$};
\node [below right of=N11] (N21) {$\bullet$};
\node [above right of=N21] (N12) {$\bullet$};
\node [below right of=N12] (N22) {$\bullet$};
\node [above right of=N22] (N13) {$\bullet$};
\node [below right of=N13] (N23) {$\bullet$};
\node [above right of=N23] (N14) {$\bullet$};
\node [below right of=N21] (N31) {$\bullet$};
\node [below right of=N22] (N32) {$\bullet$};
\node [below right of=N31] (N41) {$\bullet$};
\path[-,thick,every node/.style={font=\sffamily\small}]
(N11) edge node [below left] {} (N21)
(N12) edge node [below left] {} (N22)
%
(N21) edge node [below left] {} (N31)
(N22) edge node [below left] {} (N32)
%
;
\end{tikzpicture} 
%
\begin{tikzpicture}[node distance=0.3cm]
\node (N11) at (0,0) {$\bullet$};
\node [below right of=N11] (N21) {$\bullet$};
\node [above right of=N21] (N12) {$\bullet$};
\node [below right of=N12] (N22) {$\bullet$};
\node [above right of=N22] (N13) {$\bullet$};
\node [below right of=N13] (N23) {$\bullet$};
\node [above right of=N23] (N14) {$\bullet$};
\node [below right of=N21] (N31) {$\bullet$};
\node [below right of=N22] (N32) {$\bullet$};
\node [below right of=N31] (N41) {$\bullet$};
\path[-,thick,every node/.style={font=\sffamily\small}]
(N11) edge node [below left] {} (N21)
(N12) edge node [below left] {} (N22)
%
(N21) edge node [below left] {} (N31)
%
(N31) edge node [below left] {} (N41)
;
\end{tikzpicture} 

%
\begin{tikzpicture}[node distance=0.3cm]
\node (N11) at (0,0) {$\bullet$};
\node [below right of=N11] (N21) {$\bullet$};
\node [above right of=N21] (N12) {$\bullet$};
\node [below right of=N12] (N22) {$\bullet$};
\node [above right of=N22] (N13) {$\bullet$};
\node [below right of=N13] (N23) {$\bullet$};
\node [above right of=N23] (N14) {$\bullet$};
\node [below right of=N21] (N31) {$\bullet$};
\node [below right of=N22] (N32) {$\bullet$};
\node [below right of=N31] (N41) {$\bullet$};
\path[-,thick,every node/.style={font=\sffamily\small}]
(N11) edge node [below left] {} (N21)
(N12) edge node [below left] {} (N22)
%
(N21) edge node [below left] {} (N31)
%
;
\end{tikzpicture} 
%
\begin{tikzpicture}[node distance=0.3cm]
\node (N11) at (0,0) {$\bullet$};
\node [below right of=N11] (N21) {$\bullet$};
\node [above right of=N21] (N12) {$\bullet$};
\node [below right of=N12] (N22) {$\bullet$};
\node [above right of=N22] (N13) {$\bullet$};
\node [below right of=N13] (N23) {$\bullet$};
\node [above right of=N23] (N14) {$\bullet$};
\node [below right of=N21] (N31) {$\bullet$};
\node [below right of=N22] (N32) {$\bullet$};
\node [below right of=N31] (N41) {$\bullet$};
\path[-,thick,every node/.style={font=\sffamily\small}]
(N11) edge node [below left] {} (N21)
(N12) edge node [below left] {} (N22)
%
%
;
\end{tikzpicture} 
%
\begin{tikzpicture}[node distance=0.3cm]
\node (N11) at (0,0) {$\bullet$};
\node [below right of=N11] (N21) {$\bullet$};
\node [above right of=N21] (N12) {$\bullet$};
\node [below right of=N12] (N22) {$\bullet$};
\node [above right of=N22] (N13) {$\bullet$};
\node [below right of=N13] (N23) {$\bullet$};
\node [above right of=N23] (N14) {$\bullet$};
\node [below right of=N21] (N31) {$\bullet$};
\node [below right of=N22] (N32) {$\bullet$};
\node [below right of=N31] (N41) {$\bullet$};
\path[-,thick,every node/.style={font=\sffamily\small}]
(N11) edge node [below left] {} (N21)
%
(N21) edge node [below left] {} (N31)
%
(N31) edge node [below left] {} (N41)
;
\end{tikzpicture} 
\caption{The $11$ Kogan faces for $n=4$ whose type avoids the pattern $132$.}\label{fig:kempf}
\end{figure}
\end{example}

\subsection{Key polynomials}

One possible generalization of Schur polynomials is the so called \emph{Demazure characters},
also known as \emph{key polynomials}. The latter name was introduced by V.~Reiner and M.~Shimozono in \cite{ReinerShimozono1995},
and refers to the combinatorial model for Demazure characters introduced in \cite{Lascoux1990Keys},
where they use the so called \emph{key tableaux}.
In order to define key polynomials, we need some preliminary terminology.

Whenever $s_i \in \symS_n$ is a simple transposition, with $i \in \{1,2,\dotsc,n-1\}$,
we let $s_i$ act on $\setC[z_1,\dotsc,z_n]$ by permuting the indices:
\[
 s_i \circ f(z_1,\dotsc,z_n) = f(z_1,\dotsc,z_{i-1},z_{i+1},z_i,z_{i+2},\dotsc,z_n).
\]
Define the \emph{divided difference operator} $\partial_i$ as
\[
 \partial_i(f) = \frac{f - s_i(f)}{z_i - z_{i+1}}.
\]
Note that $\partial_i(f)$ is indeed a polynomial, and one can check that  $\partial_i(f)$ is symmetric in
the variables $z_i$ and $z_{i+1}$.
For example,
\begin{align*}
 \partial_2 (z^2_1 z_2^5 z^3_3 z_4 )  &= \frac{z^2_1 z_2^5 z^3_3 z_4 - z^2_1 z_2^3 z^5_3 z_4}{z_2-z_3} \\
 &= z^2_1z_4 \frac{z_2^5 z^3_3 -z_2^3 z^5_3}{z_2-z_3} \\
 &= z^2_1 z_2^3 z_3^3 z_4  \frac{z_2^2 -z^2_3}{z_2-z_3} \\
 &= z^2_1 z_2^3 z_3^3 z_4 (z_2 + z_3).
\end{align*}
We now define the operators $\pi_i(f) \coloneqq \partial_i(z_i f)$ for $i=1,\dotsc,n-1$
whenever $f \in \setC[z_1,\dotsc,z_n]$.
It is straightforward to verify the following properties of the $\pi_i$'s:
\begin{itemize}
 \item $\pi_i$ preserves the degree,
 \item $\pi^2_i = \pi_i$ for all $i$,
 \item $\pi_i \pi_j = \pi_j \pi_i$ whenever $|i-j|>2$,
 \item $\pi_i \pi_{i+1} \pi_i = \pi_{i+1} \pi_i  \pi_{i+1}$ for all $i$.
\end{itemize}
The last two properties allow us to make the following definition:
Let $\sigma = s_{i_1} s_{i_2} \dotsc s_{i_\ell}$ be a \emph{reduced word} of a permutation $\sigma \in \symS_n$.
Then let 
\[
 \pi_\sigma \coloneqq \pi_{i_1} \circ \pi_{i_2} \circ \dotsb \circ \pi_{i_\ell}.
\]
The action of $\pi_\sigma$ is independent of the choice of reduced word, since we have 
the relations above.

We are now ready to define the key polynomials.
Let $\lambda$ be a partition with at most $n$ parts, and let $\sigma \in \symS_n$ be a permutation.
The \emph{key polynomial} $\key_{\lambda,\sigma}(z)$ is defined as
\begin{equation}\label{eq:keyOperatorDef}
 \key_{\lambda,\sigma}(z) \coloneqq \pi_\sigma\left( z_1^{\lambda_1} \dotsm z_n^{\lambda_n} \right).
\end{equation}

\begin{example}\label{ex:operatorKey}
Let $\lambda = (2,1,0,0)$ and $\sigma = [2,4,3,1] \in \symS_4$ in one-line notation.
The permutation can be expressed as a reduced word as 
$\sigma = s_2 s_3 s_2 s_1$.
We compute the key polynomial as follows:
\begin{align*}
 \key_{\lambda,\sigma}(z) &= \pi_2 \pi_3 \pi_2 \pi_1 (z_1^2 z_2) = 
 \pi_2 \pi_3 \pi_2 \partial_1 (z_1^3 z_2 ) = \pi_2 \pi_3 \pi_2 \left( \frac{z_1^3z_2 - z_1 z^3_2}{z_1-z_2}\right) \\
 &= \pi_2 \pi_3 \pi_2  (z_1^2z_2 + z_1 z_2^2).
\end{align*}
We continue the calculation by applying $\pi_2$ and get
\begin{align*}
\key_{\lambda,\sigma}(z) &=
 \pi_2 \pi_3 (z_2 z_1^2+z_3 z_1^2+z_2^2 z_1+z_3^2 z_1+z_2 z_3 z_1).
 \end{align*}
 Applying $\pi_2 \pi_3$ then finally gives
\begin{equation}\label{eq:keyExample}
\begin{split}
\key_{\lambda,\sigma}(z) &=
z_1^2 z_2  + 
z_1^2 z_3  + 
z_1^2 z_4  + 
z_1 z_2^2  + 
z_1 z_3^2  \\
&\phantom{=}+
z_1 z_4^2  +
z_1 z_2 z_3  + 
z_1 z_2 z_4  + 
z_1 z_3 z_4 .
 \end{split}
\end{equation}
\end{example}
In general, some monomials may appear multiple times.

\bigskip 

In \cite{Kiritchenko2010}, the following formula for key polynomials using Kogan faces was proved.
This generalizes an earlier result by A.~Postnikov and R.~Stanley \cite{PostnikovStanley2008},
who considered the case covered in \cref{prop:kempfCase}.
\begin{proposition}\label{prop:keyFromFaces}
Let $\GT(\lambda,\sigma)$ be defined as the polytopal complex
\[
 \GT(\lambda,\sigma) \coloneqq  \bigcup_{\substack{\mathcal{F} \in \GT(\lambda) \\ \type(\mathcal{F}) = w_0\sigma}} \mathcal{F}.
\]
That is, $\GT(\lambda,\sigma)$ is the union of all reduced Kogan faces of type $w_0\sigma$
in the polytope $\GT(\lambda)$.
The key polynomial $\key_{\lambda,\sigma}(z)$ can be computed as
\begin{align}\label{eq:keyGTFormula}
 \key_{\lambda,\sigma}(z_1,\dotsc,z_n) = \sum_{G \in \GT(\lambda,\sigma) \cap \setZ^{\tfrac{n(n+1)}{2}} } z_1^{w_1(G)} \dotsm z_n^{w_n(G)}
\end{align}
where we use the same weight for integral Gelfand--Tsetlin patterns as in \eqref{eq:weight}.
\end{proposition}
As an immediate corollary, it is clear that
\[
 \schurS_\lambda(z_1,\dotsc,z_n) = \key_{\lambda,w_0}(z_1,\dotsc,z_n).
\]

We now recalculate the key polynomial in \cref{ex:operatorKey}
using \eqref{eq:keyGTFormula}.
\begin{example}\label{ex:gtKey}
Let $\lambda = (2,1,0,0)$ and $\sigma = [2,4,3,1] \in \symS_4$.
We have that $\omega_0 \sigma = [1,3,4,2]$,
and we have that 
\[
 [1,3,4,2] = s_3 s_2.
\]
There are no other reduced words that give rise to the same permutation.
However, there are three reduced Kogan faces that give rise to this particular reduced word (and are hence of type $[1,3,4,2]$):
\begin{equation}\label{eq:gtFaces}
\begin{tikzpicture}[node distance = 0.8cm]
\node (N11) at (0,0) {$\bullet$};
\node [below right of=N11] (N21) {$\bullet$};
\node [above right of=N21] (N12) {$\bullet$};
\node [below right of=N12] (N22) {$\bullet$};
\node [above right of=N22] (N13) {$\bullet$};
\node [below right of=N13] (N23) {$\bullet$};
\node [above right of=N23] (N14) {$\bullet$};
\node [below right of=N21] (N31) {$\bullet$};
\node [below right of=N22] (N32) {$\bullet$};
\node [below right of=N31] (N41) {$\bullet$};
\node [below of=N41] (L1) {$(A)$};

\path[-,thick,every node/.style={font=\sffamily\small}]
(N12) edge node [below left] {$s_2$} (N22)
(N22) edge node [below left] {$s_3$} (N32);
\end{tikzpicture}
\quad
 \begin{tikzpicture}[node distance = 0.8cm]
\node (N11) at (0,0) {$\bullet$};
\node [below right of=N11] (N21) {$\bullet$};
\node [above right of=N21] (N12) {$\bullet$};
\node [below right of=N12] (N22) {$\bullet$};
\node [above right of=N22] (N13) {$\bullet$};
\node [below right of=N13] (N23) {$\bullet$};
\node [above right of=N23] (N14) {$\bullet$};
\node [below right of=N21] (N31) {$\bullet$};
\node [below right of=N22] (N32) {$\bullet$};
\node [below right of=N31] (N41) {$\bullet$};
\node [below of=N41] (L1) {$(B)$};
\path[-,thick,every node/.style={font=\sffamily\small}]
(N31) edge node [below left] {$s_3$} (N41)
(N12) edge node [below left] {$s_2$} (N22);
\end{tikzpicture}
\quad
 \begin{tikzpicture}[node distance = 0.8cm]
\node (N11) at (0,0) {$\bullet$};
\node [below right of=N11] (N21) {$\bullet$};
\node [above right of=N21] (N12) {$\bullet$};
\node [below right of=N12] (N22) {$\bullet$};
\node [above right of=N22] (N13) {$\bullet$};
\node [below right of=N13] (N23) {$\bullet$};
\node [above right of=N23] (N14) {$\bullet$};
\node [below right of=N21] (N31) {$\bullet$};
\node [below right of=N22] (N32) {$\bullet$};
\node [below right of=N31] (N41) {$\bullet$};
\node [below of=N41] (L1) {$(C)$};
\path[-,thick,every node/.style={font=\sffamily\small}]
 (N21) edge node [below left] {$s_2$} (N31)
(N31) edge node [below left] {$s_3$} (N41);
\end{tikzpicture}
\end{equation}
We expect nine lattice points in the union of these faces,
as there are nine monomials in \eqref{eq:keyExample}.
These lattice points are given by the following Gelfand--Tsetlin patterns:
\[
\begin{tikzpicture}[scale=0.4]
\node at (3,0) {$2$};
\node at (5,0) {$1$};
\node at (7,0) {$0$};
\node at (9,0) {$0$};
\node at (4,-1) {$1$};
\node at (6,-1) {$0$};
\node at (8,-1) {$0$};
\node at (5,-2) {$1$};
\node at (7,-2) {$0$};
\node at (6,-3) {$1$};
\node at (6,-4) {$C$};
\node at (6,-5) {$z_1z_4^2$};
\end{tikzpicture}
\qquad
\begin{tikzpicture}[scale=0.4]
\node at (3,0) {$2$};
\node at (5,0) {$1$};
\node at (7,0) {$0$};
\node at (9,0) {$0$};
\node at (4,-1) {$1$};
\node at (6,-1) {$1$};
\node at (8,-1) {$0$};
\node at (5,-2) {$1$};
\node at (7,-2) {$0$};
\node at (6,-3) {$1$};
\node at (6,-4) {$BC$};
\node at (6,-5) {$z_1z_3z_4$};
\end{tikzpicture}
\qquad
\begin{tikzpicture}[scale=0.4]
\node at (3,0) {$2$};
\node at (5,0) {$1$};
\node at (7,0) {$0$};
\node at (9,0) {$0$};
\node at (4,-1) {$1$};
\node at (6,-1) {$1$};
\node at (8,-1) {$0$};
\node at (5,-2) {$1$};
\node at (7,-2) {$1$};
\node at (6,-3) {$1$};
\node at (6,-4) {$ABC$};
\node at (6,-5) {$z_1z_2z_4$};
\end{tikzpicture}
\]
\[
\begin{tikzpicture}[scale=0.4]
\node at (3,0) {$2$};
\node at (5,0) {$1$};
\node at (7,0) {$0$};
\node at (9,0) {$0$};
\node at (4,-1) {$2$};
\node at (6,-1) {$0$};
\node at (8,-1) {$0$};
\node at (5,-2) {$2$};
\node at (7,-2) {$0$};
\node at (6,-3) {$2$};
\node at (6,-4) {$C$};
\node at (6,-5) {$z_1^2z_4$};
\end{tikzpicture}
\qquad
\begin{tikzpicture}[scale=0.4]
\node at (3,0) {$2$};
\node at (5,0) {$1$};
\node at (7,0) {$0$};
\node at (9,0) {$0$};
\node at (4,-1) {$2$};
\node at (6,-1) {$1$};
\node at (8,-1) {$0$};
\node at (5,-2) {$1$};
\node at (7,-2) {$0$};
\node at (6,-3) {$1$};
\node at (6,-4) {$B$};
\node at (6,-5) {$z_1z_3^2$};
\end{tikzpicture}
\qquad
\begin{tikzpicture}[scale=0.4]
\node at (3,0) {$2$};
\node at (5,0) {$1$};
\node at (7,0) {$0$};
\node at (9,0) {$0$};
\node at (4,-1) {$2$};
\node at (6,-1) {$1$};
\node at (8,-1) {$0$};
\node at (5,-2) {$1$};
\node at (7,-2) {$1$};
\node at (6,-3) {$1$};
\node at (6,-4) {$AB$};
\node at (6,-5) {$z_1z_2z_3$};
\end{tikzpicture}
\]
\[
\begin{tikzpicture}[scale=0.4]
\node at (3,0) {$2$};
\node at (5,0) {$1$};
\node at (7,0) {$0$};
\node at (9,0) {$0$};
\node at (4,-1) {$2$};
\node at (6,-1) {$1$};
\node at (8,-1) {$0$};
\node at (5,-2) {$2$};
\node at (7,-2) {$0$};
\node at (6,-3) {$2$};
\node at (6,-4) {$BC$};
\node at (6,-5) {$z_1^2z_3$};
\end{tikzpicture}
\qquad
\begin{tikzpicture}[scale=0.4]
\node at (3,0) {$2$};
\node at (5,0) {$1$};
\node at (7,0) {$0$};
\node at (9,0) {$0$};
\node at (4,-1) {$2$};
\node at (6,-1) {$1$};
\node at (8,-1) {$0$};
\node at (5,-2) {$2$};
\node at (7,-2) {$1$};
\node at (6,-3) {$1$};
\node at (6,-4) {$A$};
\node at (6,-5) {$z_1z_2^2$};
\end{tikzpicture}
\qquad
\begin{tikzpicture}[scale=0.4]
\node at (3,0) {$2$};
\node at (5,0) {$1$};
\node at (7,0) {$0$};
\node at (9,0) {$0$};
\node at (4,-1) {$2$};
\node at (6,-1) {$1$};
\node at (8,-1) {$0$};
\node at (5,-2) {$2$};
\node at (7,-2) {$1$};
\node at (6,-3) {$2$};
\node at (6,-4) {$AC$};
\node at (6,-5) {$z_1^2z_2$};
\end{tikzpicture}
\]
The letters below each pattern indicate which reduced Kogan faces in \eqref{eq:gtFaces}
the pattern is a member of and the monomials represent $z^{w(G)}$ as defined in \eqref{eq:weight}.
\end{example}

\section{Ehrhart polynomials}\label{sec:ehrhart}

From \eqref{eq:weylDimFormula}, it follows that the Ehrhart polynomial of $\GT(\lambda)$ is given by
\begin{equation} \label{eq:schurEhrhart}
\ehr(\GT(\lambda),k) = \prod_{1 \leq i < j \leq n} \frac{k(\lambda_i - \lambda_j) + j-i}{j-i} 
\end{equation}
and it is clear that all coefficients of $k$ are non-negative.

\medskip 
For a skew shape $\lambda/\mu$, the lattice point enumerator is indeed a polynomial  $\ehr(\GT(\lambda/\mu),k) = \schurS_{k\lambda/k\mu}(\underbrace{1,1,\dotsc,1}_n)$. Table 1 below shows a sample computation of such polynomials. Computer evidence suggests the following conjecture:

\begin{conjecture}
Let $\lambda/\mu$ be a skew shape. Then the polynomial $\ehr(\GT(\lambda/\mu),k)$ has non-negative coefficients.
\end{conjecture}
We do not expect a closed-form formula for $\ehr(\GT(\lambda/\mu),k)$ --- 
the volume of $\GT(\lambda/\mu)$ is related to the 
number of skew standard Young tableaux for which there are 
no known closed formulas in general, see \cite{MoralesPakPanova2018} and subsequent papers.

\begin{table}[!ht]
\centering
{\begin{tabular}{@{\extracolsep{4pt}}cl}
$\lambda/\mu$ & $\schurS_{k\lambda/k\mu}(1,1,1)$ \\
\toprule 
\vspace{0.3cm}
$221$ &  $\frac{1}{2} \left(k^2+3 k+2\right)$ \\
\vspace{0.3cm}
$\substack{221/1 \\ 221/21 \\ 221/11} $  & $\frac{1}{4} \left(k^4+6 k^3+13 k^2+12 k+4\right)$ \\
\vspace{0.3cm}
$\substack{221/2 \\ 321 }$  & $k^3+ 3 k^2 + 3 k +1$ \\
\vspace{0.3cm}
 $\substack{321/1 \\ 321/2 }$ & $\frac{1}{2} (k^5+6 k^4+14 k^3+16 k^2+9 k+2)$ \\
\vspace{0.3cm}
 $321/21$ & $\frac18  ( k^6+ 9 k^5+ 33 k^4 + 63 k^3+ 66 k^2+ 36 k +8)$ \\
\bottomrule
\end{tabular}}
\caption{Ehrhart polynomials for $\GT(\lambda/\mu)$.}
\end{table}

\subsection{Kostka coefficients}\label{ssec:kostka}

Kostka coefficients are important numbers appearing in several branches of 
mathematics such as algebraic combinatorics, symmetric functions, representation 
theory and algebraic geometry among others. 
From a representation theoretical point of view, they are defined as 
the dimension of the weight subspace of the irreducible representation of 
the Lie algebra $\mathfrak{gl}(\setC)$.
From a combinatorial point of view, they enumerate the number of 
semi-standard Young tableaux of fixed shape and weight.

Given two partitions $\lambda$ and $\mu$ of $n$, Kostka coefficients arise in the 
expansion of the Schur polynomial as a linear combination of monomial symmetric polynomials
\begin{equation}
\schurS_\lambda(z_1,\dotsc,z_n)=\sum_{\mu} K_{\lambda \mu} m_\mu(z_1,\dots,z_n).
\end{equation}

More generally, Kostka coefficients can be viewed as a special case of so called Littlewood--Richardson coefficients
 $K_{\lambda\mu}=c_{\sigma\lambda}^\tau$, where $\sigma$ and $\tau$ are defined 
 in terms of the partition $\mu$. 
 For more material about the latter coefficients and some standard results, 
 the reader is referred to the book \cite{Fulton91representation}.

The partitions $\lambda$ and $\mu$ are usually represented as integer vectors, 
so it makes sense to talk about scaling these vectors by an integer factor $k$. 
This operation gives rise to what is called ``stretched'' Kostka coefficients 
$K_{\lambda\mu}(k)\coloneqq K_{k\lambda, k\mu}$.

Surprisingly at first, $K_{\lambda\mu}(k)$ turns out to be a polynomial
function in $k$ with rational coefficients depending on $\lambda$ and $\mu$ --- 
a fact initially shown by A.~Kirillov and N.~Reshetikhin \cite{Kirillov88thebethe}.

Expressing Kostka coefficients using the Gelfand--Tsetlin polytopes 
from \cref{ssec:gtPolytope} provides a 
natural geometric interpretation of $K_{\lambda\mu}(k)$. 
To each Kostka coefficient $K_{\lambda\mu}$, 
we have the corresponding Gelfand--Tsetlin 
polytope $\GT(\lambda,\mu)$ in $\mathbb{R}^\frac{n(n+1)}{2}$ such that
\begin{equation}
K_{\lambda\mu}=|\GT(\lambda,\mu)\cap \mathbb{Z}^\frac{n(n+1)}{2}| \text{ and }
 K_{\lambda\mu}(k) = |k\cdot \GT(\lambda,\mu)\cap \mathbb{Z}^\frac{n(n+1)}{2}|.
\end{equation}
The vertices of $\GT(\lambda,\mu)$ have rational coordinates in general $(n\geq5)$. 
From E.~Ehrhart's fundamental work, it is well-known that $K_{\lambda\mu}(k)$ must be a quasipolynomial. 
This means that there exist an integer $M$ and polynomials $g_0,g_1,\dots,g_{M-1}$ 
such that $K_{\lambda\mu}(k)=g_i(k)$ whenever $k\equiv i\ (\textrm{mod}\ M )$ --- see details in \cite{StanleyEC1}. 
The ``surprising'' fact here is that the function $K_{\lambda\mu}(k)$ is indeed a polynomial, 
which exhibits an example of \emph{period collapse} in rational polytopes.

In their paper \cite{King04stretched}, the authors conjectured a formula for 
the degree of $K_{\lambda\mu}(k)$ which was later proved by T.~McAllister \cite{Mcallister2008}. 
In that same paper, they also conjectured that the 
coefficients of the polynomial $K_{\lambda\mu}(k)$ are non-negative. 
To the best of our knowledge, this conjecture remains open and provides an 
instance of a well-known phenomenon called Ehrhart positivity. 
For a recent survey about this topic, we recommend the article 
by F.~Liu \cite{Liu2017}.

\medskip

The polytopes $\GT(\lambda/ \mu, w)$ also exhibit a period collapse, and have a polynomial Ehrhart funtion.
This follows from the work of E.~Rassart in \cite{Rassart2004}.
Integrality and the integer decomposition property of the family $\GT(\lambda/ \mu, w)$ is studied
in the work of P.~Alexandersson \cite{Alexandersson2016GTPoly}, where several questions are left unanswered.
The \emph{skew Kostka coefficients} $K_{\lambda/\mu,\nu}$ are defined analogously to the usual Kostka coefficients, 
via the identity
\begin{equation}
\schurS_{\lambda/\mu}(z_1,\dotsc,z_n)=\sum_{\nu} K_{\lambda/\mu, \nu} m_\nu(z_1,\dots,z_n).
\end{equation}

Scaling by an integer factor $k$ gives rise to the polynomial $\ehr(\GT(\lambda/\mu,\nu),k) = K_{k\lambda/k\mu,k\nu} \in \setQ[k]$. The conjecture by R.C.~King et al. seems to extend to the skew case:
\begin{conjecture}
Let $\lambda/\mu$ be a skew shape and $\nu$ be a weight. Then the polynomial $\ehr(\GT(\lambda/\mu,\nu),k)$ has non-negative coefficients.
\end{conjecture}


\section{A conjecture on key polynomials}\label{sec:conjectures}

Given the positivity phenomena above, it is reasonable to ask if 
we have positive Ehrhart coefficients for the polytopal complexes in \cref{prop:keyFromFaces}. 
In other words, for a partition $\lambda$ and permutation $\sigma \in \symS_n$,
do
\[
\ehr(\GT(\lambda,\sigma),k) = \key_{k\lambda,\sigma}(\underbrace{1,1,\dotsc,1}_n) \in \setQ[k]
\]
have non-negative coefficients?

By using the divided difference operators in \cref{eq:keyOperatorDef}, it is possible to compute the Ehrhart polynomials
for small values of $n$.
Note first that for all $\lambda$, $\sigma$ and $m \geq 0$, we have that
\begin{equation}\label{eq:factorKey}
\key_{m+\lambda,\sigma}(z_1,\dotsc,z_n) = (z_1z_2\dotsm z_n)^m \key_{\lambda,\sigma}(z_1,\dotsc,z_n),
\end{equation}
where $m+\lambda \coloneqq (m+\lambda_1,m+\lambda_2,\dotsc,m+\lambda_n)$. 
It follows that
\[
\key_{m+\lambda,\sigma}(1,\dotsc,1)=\key_{\lambda,\sigma}(1,\dotsc,1) \text{ and }
\ehr(\GT(m+\lambda,\sigma),k) = \ehr(\GT(\lambda,\sigma),k).
\]

\begin{example}
Let $\lambda = (a+b,a,0)$, with $a,b\geq 0$, since we can assume that the last part is equal to zero, due to \eqref{eq:factorKey}.
We have the following table for $\key_{\lambda,\sigma}(z_1,z_2,z_3)$, for different values of $\sigma \in \symS_3$:
\begin{table}[!ht]
\centering
{\begin{tabular}{@{\extracolsep{4pt}}cl}
$\sigma$ & $\key_{\lambda,\sigma}(z_1,z_2,z_3)$ \\
\toprule 
\vspace{0.3cm}
$[1,2,3]$ & $z^{a+b}_1 z^{a}_2$ \\
\vspace{0.3cm}
$[2,1,3]$ & $\frac{z_1^a z_2^a \left(z_1^{b+1}-z_2^{b+1}\right)}{z_1-z_2} $ \\
\vspace{0.3cm}
$[1,3,2]$ & $\frac{\left(z_2^{a+1}-z_3^{a+1}\right) z_1^{a+b}}{z_2-z_3} $ \\
\vspace{0.3cm}
$[3,1,2]$ & $
\frac{
z_1^a \left(\left(z_1-z_3\right) z_2^{a+1} \left(z_1^{b+1}-z_2^{b+1}\right)-\left(z_1-z_2\right) z_3^{a+1} \left(z_1^{b+1}-z_3^{b+1}\right)\right)
}{\left(z_1-z_2\right) \left(z_1-z_3\right) \left(z_2-z_3\right)}
$ \\
\vspace{0.3cm}
$[2,3,1]$ & $
\frac{\left(z_3-z_1\right) \left(z_3^{a+1}-z_2^{a+1}\right) z_1^{a+b+1}+\left(z_2-z_3\right) \left(z_3^{a+1}-z_1^{a+1}\right) z_2^{a+b+1}}{\left(z_1-z_2\right) \left(z_1-z_3\right) \left(z_2-z_3\right)}
$ \\
\vspace{0.3cm}
$[3,2,1]$ & $
\frac{z_1^{a+1} \left(z_3^{a+b+2}-z_2^{a+b+2}\right)+\left(z_2^{a+1}-z_3^{a+1}\right) z_1^{a+b+2}+z_2^{a+1} z_3^{a+1} \left(z_2^{b+1}-z_3^{b+1}\right)}{\left(z_1-z_2\right) \left(z_1-z_3\right) \left(z_2-z_3\right)}\label{tab:keypolys}
$ \\
\bottomrule
\end{tabular}}
\caption{Key polynomials computed for all $\sigma \in \symS_3$.}
\end{table}
\medskip 

By taking the limit $z_i \to 1$ in \cref{tab:keypolys}, and multiplying $a$ and $b$ with $k$,
we get the Ehrhart polynomials $\ehr(\GT(\lambda,\sigma),k)$, which we present in \cref{tab:keyEhrhart}.
\begin{table}[!ht]
\centering
{\begin{tabular}{@{\extracolsep{60pt}}cl}
$\sigma$ & $\key_{k\lambda,\sigma}(1,1,1) = \ehr(\GT(\lambda,\sigma),k)$ \\
\toprule 
\vspace{0.3cm}
$[1,2,3]$ & $1$ \\
\vspace{0.3cm}
$[2,1,3]$ & $1+bk $ \\
\vspace{0.3cm}
$[1,3,2]$ & $1+ak $ \\
\vspace{0.3cm}
$[3,1,2]$ & $\frac{1}{2}(bk+1) (2ak + bk+2)$ \\
\vspace{0.3cm}
$[2,3,1]$ & $ \frac{1}{2} (a k+1) (2bk + ak+2) $ \\
\vspace{0.3cm}
$[3,2,1]$ & $\frac12 (ak+1) (b k+1) (2+a k+b k)$ \\
\bottomrule
\end{tabular}}
\caption{Values of $\key_{k\lambda,\sigma}(z_1,z_2,z_3)$ after taking the limit $z_i\to 1$. Note that
$a$ and $b$ have been multiplied by $k$.}
\label{tab:keyEhrhart}
\end{table}
The last polynomial for $\sigma=[3,2,1]$ agrees with \cref{eq:schurEhrhart} with $\lambda_1=a+b$, $\lambda_2 = a$ and $\lambda_3=0$, since we get
\[
\key_{k \lambda,[3,2,1]}(1,1,1) = 
  \frac{k(a+b - a) + (2-1)}{2-1} 
  \frac{k(a - 0) + (3-2)}{3-2} 
  \frac{k(a+b - 0) + (3-1)}{3-1}.
\]
\end{example}

We are then fairly confident in the following conjecture:
\begin{conjecture}\label{conj:keyEhrhart}
 The Ehrhart polynomial $\ehr(\GT(\lambda,\sigma),k) \in \setQ[k]$ has only non-negative coefficients
 for all $\lambda$, $\sigma$.
Furthermore, in the new variables variables $\lambda_i = a_1+a_2+\dotsb + a_i$, we have
\begin{equation}
 \ehr(\GT(\lambda,\sigma),k) \in \setQ[k,a_1,a_2,\dotsc,a_n]
\end{equation}
with non-negative coefficients.
\end{conjecture}
Above, we verified this for $\sigma \in \symS_3$, and we have also verified this for $\sigma \in \symS_4$.
The symbolic computations done in Mathematica become tedious for larger values of $n$.

\subsection{A determinant formula}\label{ssec:determinant}

An interesting special case was considered in \cite{PostnikovStanley2008}, where $\sigma$ is $231$-avoiding and 
$\ehr(\GT(\lambda,\sigma),k)$ reduces to the Ehrhart polynomial of a single 
reduced Kogan face of $\GT(\lambda)$, with type $w_0 \sigma$.
The number of $231$-avoiding permutations in $\symS_n$ is given by the Catalan numbers, $\frac{1}{n+1}\binom{2n}{n}$.
For an excellent survey on Catalan numbers, see the book by R.~Stanley \cite{StanleyCatalan}.

From \cite[Corollary 14.6]{PostnikovStanley2008}, it follows that the Ehrhart polynomial
$\ehr(\GT(\lambda,\sigma),k)$ for a $231$-avoiding permutation $\sigma\in \symS_n$ can be obtained as the determinant
\begin{equation}\label{eq:determinant}
 \ehr(\GT(\lambda,\sigma),k) = \det\left(  \binom{k\lambda_i + b_i - i}{ b_i - j}  \right)_{1\leq i,j \leq n}
\end{equation}
where $b_1,b_2,\dotsc,b_n$ form a sequence of non-negative integers (determined by $\sigma$) satisfying 
\[
b_1 \leq b_2 \leq \dotsb \leq b_n \leq n \text{ and } b_i\geq i \text{ for } i=1,2,\dotsc,n.
\]
Such sequences are in bijection with $231$-avoiding permutations in $\symS_n$.
Conjecture~\ref{conj:keyEhrhart} thus states that for every integer partition $\lambda$,
the determinant in \eqref{eq:determinant} has non-negative coefficients (as a polynomial in $k$).

The reason \eqref{eq:determinant} exists is due to a Jacobi--Trudi identity for flagged Schur polynomials,
as key polynomials are certain flagged Schur polynomials in this special setting.

\section{Faces of Gelfand--Tsetlin polytopes}\label{sec:faces}

As a consequence of Conjecture~\ref{conj:keyEhrhart} and \cref{prop:kempfCase}, it would follow
that certain faces of Gelfand--Tsetlin polytopes have non-negative coefficients in the Ehrhart polynomial.
It is therefore appropriate to investigate Ehrhart polynomials of general faces of $\GT(\lambda)$.

However, in this generality there are certain faces where negative coefficients appear.
The following example was constructed in \cite{Alexandersson16Counterexamples}:
\begin{example}
Consider the following face $\mathcal{F}_\ell$ of $\GT(\lambda)$ where $\lambda = (1^\ell, 0^{\ell+1})$:
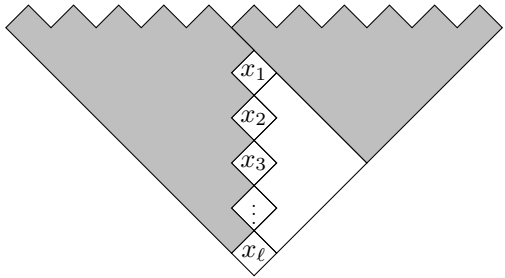
\[
\begin{tikzpicture}[scale=0.3]
\draw[black] (12,-2)--(13,-3)--(14,-2)--(13,-1)--(12,-2)--cycle;
\draw[black] (13,-9)--(14,-10)--(15,-9)--(16,-8)--(17,-7)--(18,-6)--(17,-5)--(16,-4)--(15,-3)--(14,-2)--(13,-3)--(14,-4)--(13,-5)--(14,-6)--(13,-7)--(14,-8)--(13,-9)--cycle;
\draw[black] (12,-4)--(13,-5)--(14,-4)--(13,-3)--(12,-4)--cycle;
\draw[black] (12,-6)--(13,-7)--(14,-6)--(13,-5)--(12,-6)--cycle;
\draw[black] (12,-8)--(13,-9)--(14,-8)--(13,-7)--(12,-8)--cycle;
\draw[black] (12,-10)--(13,-11)--(14,-10)--(13,-9)--(12,-10)--cycle;
\draw[color=black,fill=lightgray] (18,-6)--(19,-5)--(20,-4)--(21,-3)--(22,-2)--(23,-1)--(24,0)--(23,1)--(22,0)--(21,1)--(20,0)--(19,1)--(18,0)--(17,1)--(16,0)--(15,1)--(14,0)--(13,1)--(12,0)--(13,-1)--(14,-2)--(15,-3)--(16,-4)--(17,-5)--(18,-6)--cycle;
\draw[color=black,fill=lightgray] (12,-10)--(13,-9)--(12,-8)--(13,-7)--(12,-6)--(13,-5)--(12,-4)--(13,-3)--(12,-2)--(13,-1)--(12,0)--(11,1)--(10,0)--(9,1)--(8,0)--(7,1)--(6,0)--(5,1)--(4,0)--(3,1)--(2,0)--(3,-1)--(4,-2)--(5,-3)--(6,-4)--(7,-5)--(8,-6)--(9,-7)--(10,-8)--(11,-9)--(12,-10)--cycle;
\node at (3,0) {$ $};
\node at (5,0) {$ $};
\node at (7,0) {$ $};
\node at (9,0) {$ $};
\node at (11,0) {$ $};
\node at (13,0) {$ $};
\node at (15,0) {$ $};
\node at (17,0) {$ $};
\node at (19,0) {$ $};
\node at (21,0) {$ $};
\node at (23,0) {$ $};
\node at (4,-1) {$ $};
\node at (6,-1) {$ $};
\node at (8,-1) {$ $};
\node at (10,-1) {$ $};
\node at (12,-1) {$ $};
\node at (14,-1) {$ $};
\node at (16,-1) {$ $};
\node at (18,-1) {$ $};
\node at (20,-1) {$ $};
\node at (22,-1) {$ $};
\node at (5,-2) {$ $};
\node at (7,-2) {$ $};
\node at (9,-2) {$ $};
\node at (11,-2) {$ $};
\node at (13,-2) {$x_1$};
\node at (15,-2) {$ $};
\node at (17,-2) {$ $};
\node at (19,-2) {$ $};
\node at (21,-2) {$ $};
\node at (6,-3) {$ $};
\node at (8,-3) {$ $};
\node at (10,-3) {$ $};
\node at (12,-3) {$ $};
\node at (14,-3) {$ $};
\node at (16,-3) {$ $};
\node at (18,-3) {$ $};
\node at (20,-3) {$ $};
\node at (7,-4) {$ $};
\node at (9,-4) {$ $};
\node at (11,-4) {$ $};
\node at (13,-4) {$x_2$};
\node at (15,-4) {$ $};
\node at (17,-4) {$ $};
\node at (19,-4) {$ $};
\node at (8,-5) {$ $};
\node at (10,-5) {$ $};
\node at (12,-5) {$ $};
\node at (14,-5) {$ $};
\node at (16,-5) {$ $};
\node at (18,-5) {$ $};
\node at (9,-6) {$ $};
\node at (11,-6) {$ $};
\node at (13,-6) {$x_3$};
\node at (15,-6) {$ $};
\node at (17,-6) {$ $};
\node at (10,-7) {$ $};
\node at (12,-7) {$ $};
\node at (14,-7) {$ $};
\node at (16,-7) {$ $};
\node at (11,-8) {$ $};
\node at (13,-8) {$ \scriptstyle{\vdots} $};
\node at (15,-8) {$ $};
\node at (12,-9) {$ $};
\node at (14,-9) {$ $};
\node at (13,-10) {$x_\ell$};
\end{tikzpicture}
\]
All entries in each respective region are forced to be equal by additional constraints, thus giving rise to a face of 
a Gelfand--Tsetlin polytope.
One can construct a bijection between lattice points in this face, and lattice points in a certain order polytope.
It is then straightforward to prove that
\[
 \ehr(\mathcal{F}_\ell,k) = \sum_{j=1}^{k+1} j^\ell \in \setQ[k]
\]
and this polynomial has a negative coefficient for $\ell = 20$.
\end{example}
For more background on negative Ehrhart coefficients of order polytopes, we refer to \cite{LiuTsuchiya2018}.

\section*{Acknowledgements}

The authors would like to thank Osaka University for their hospitality, in particular T.~Hibi and A.~Tsuchiya. The first author is funded by the \emph{Knut and Alice Wallenberg Foundation} (2013.03.07).

 \bibliographystyle{amsalpha}
 \bibliography{bibliography}

\end{document}